\documentclass[11pt, a4paper]{report}
\usepackage[colorlinks]{hyperref}
\usepackage{amsmath, amsfonts, amssymb, amsthm, overpic, geometry, sgamevar}
\DeclareMathOperator*{\argmax}{arg\,max}
\DeclareMathOperator*{\argmin}{arg\,min}
\DeclareMathOperator*{\Image}{Im}
\DeclareMathOperator*{\isom}{isom}
\DeclareMathOperator*{\inv}{inv}
\DeclareMathOperator*{\bij}{bij}
\DeclareMathOperator*{\Aut}{Aut}

\DeclareMathOperator*{\id}{id}

\newtheoremstyle{thesis}
{3pt}
{3pt}
{}
{}
{\bfseries}
{:}
{.5em}
{}

\theoremstyle{thesis}
\newtheorem{theorem}{Theorem}[section]
\newtheorem{definition}[theorem]{Definition}
\newtheorem{conjecture}[theorem]{Conjecture}
\newtheorem{example}[theorem]{Example}
\newtheorem{proposition}[theorem]{Proposition}
\newtheorem{definition-proposition}[theorem]{Definition-Proposition}
\newtheorem{properties}[theorem]{Properties}
\newtheorem{remark}[theorem]{Remark}

\begin{document}
	\begin{titlepage}
	\begin{center}
		\vspace*{3cm}
		
		{\huge\sffamily\bfseries Some Structure Properties of \\ Finite Normal-Form Games}
		
		\vspace{0.7cm}
		
		{\large by}
		
		\vspace{0.1cm}
		
		{\Large Nicholas Ham}
		
		\vspace{0.1cm}
				
		{\large BEc, BSc}
		
		\vspace{0.7cm}
		
  		{\Large Supervised by Des FitzGerald}
  
  		\vspace{1.8cm}
  		\includegraphics[scale = 1.1]{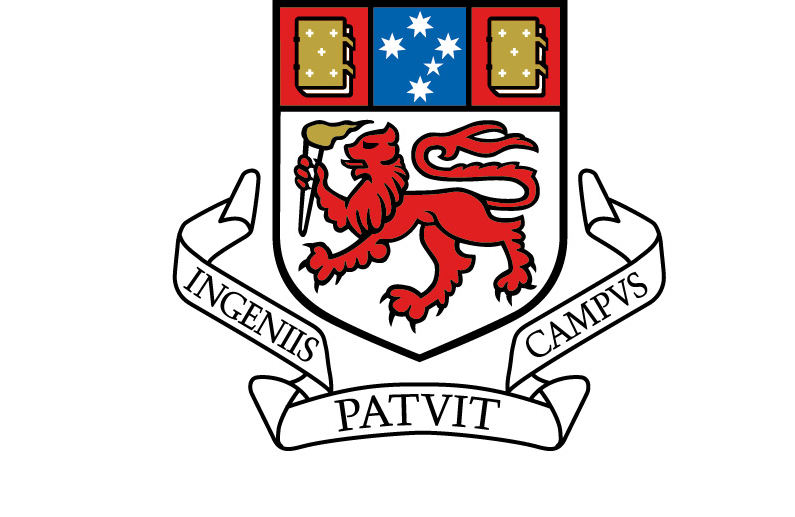}
  		
  		\vspace{1.4cm}
      
		{\large A thesis submitted in fulfilment of the requirements}
		
		\vspace{0.25cm}
		
		{\large for the Degree of Bachelor of Science with Honours}
		
		\vspace{0.7cm}
		
		{\large School of Mathematics \& Physics}
		
		\vspace{0.25cm}
		
		{\large University of Tasmania}
		
		\vspace{0.7cm}
		
		{\large 2011}
		
		\vspace{1.5cm}
	\end{center}
\end{titlepage}

    \begin{abstract}
        Game theory provides a mathematical framework for analysing strategic situations involving at least two players. Normal-form games model situations where the players simultaneously pick their moves. In this thesis we explore the strategic structure of finite normal-form games.
        
        We look at three notions of isomorphisms between games, the structural properties that they preserve and under what conditions they are met. We also look at various notions of symmetric games, under what conditions they are met, the structural properties that these notions capture, how to identify them and how to construct them. 
    \end{abstract}

    \tableofcontents
    
    \chapter{Introduction} 
Game theory provides a mathematical framework for analysing strategic situations involving at least two players. Since the foundational work of von Neumann and Morgenstern \cite{VNM} it has found applications in areas such as artificial intelligence, biology and economics. 

A normal-form game models a situation where the players simultaneously pick their strategy from a set of strategies and each player assigns a utility or payoff value to every possible combination of the strategies.

Most game theory literature is concerned with how a rational player can optimally select their strategy to maximise their utility, given that each of their opponents is attempting to do the same thing. A lot of this analysis is concerned with the classification of equilibrium concepts and proving their existence. We take an excursion away from this pursuit and explore the strategic structure of finite normal-form games.

Analysing the strategic structure of games should be of interest not only to mathematicians but also people working in other disciplines, for things like determining if a game is fair or under what conditions two situations arising from seemingly different contexts could be considered to have the same strategic structure.

In Chapter \ref{chap:background} we familiarise the reader with some notation and common mathematical terms which will be needed throughout the remainder of the thesis. We provide a proof without reference for any results that may not be common knowledge although we claim no originality to these results.

In Chapter \ref{chap:games} we give a complete definition of a game and a brief introduction to several concepts which will play a central role to the strategic structure of games we want to analyse. 

In Chapter \ref{chap:equivalence} we begin by looking at three notions of isomorphisms between games, the structural properties that they preserve and under what conditions they are met. In particular we look at the necessary conditions for the preservation of structure relating to the pure and mixed strategy spaces of two games. We finish the chapter by defining equivalence between games for each of these notions.

In Chapter \ref{chap:symmgames} we move our focus to symmetric games. We begin by reviewing various notions of a game being symmetric and introduce classifications which better capture the different possible types of symmetric structure. We look at the conditions required for these notions to be met, the structural properties that these notions capture, how to identify them, how to construct them and finally look at two famous theorems on the existence of equilibria in symmetric games.
    \chapter{Background} \label{chap:background}
    We begin by familiarising the reader with some notation and mathematical terms that will be used throughout the remainder of the thesis. We prove some of the results that may not be common knowledge, but we claim no originality to them.
        
    \section{Binary Relations}
        \begin{definition} 
            A \textit{binary relation} $R$ between two sets $A$ and $B$ is a subset of $A\times B$. We say that $a \in A$ is \textit{related} to $b \in B$ if $(a,b) \in R$ and write $a R b$. The $inverse$ of a binary relation $R$ is $\{(b,a): (a,b) \in R\} \subseteq B\times A$ denoted as $R^{-1}$.
        \end{definition}
            
        \begin{properties} 
            A binary relation $R \subseteq A\times B$ is; 
            \begin{itemize}
                \item \textit{injective} if for all $a_1, a_2 \in A$ and $b \in B$, $(a_1, b) \in R$ and $(a_2, b) \in R$ implies $a_1 = a_2$,
                \item \textit{functional} if for all $a \in A$ and $b_1, b_2 \in B$, $(a, b_1), (a, b_2) \in R$ implies $b_1 = b_2$,
                \item \textit{left-total} if for all $a \in A$ there exists $b \in B$ such that $(a,b) \in R$,  
                \item \textit{surjective} if for all $b \in B$ there exists $a \in A$ such that $(a,b) \in R$,
                \item \textit{reflexive} if $B = A$ and $(a, a) \in R$ for all $a \in A$,
                \item \textit{symmetric} if $B = A$ and for all $a_1, a_2 \in A$, $(a_1, a_2) \in R$ implies $(a_2, a_1) \in R$.
                \item \textit{transitive} if $B = A$ and $(a_1, a_2), (a_2, a_3) \in R$ implies $(a_1, a_3) \in R$, and
                \item \textit{total} if $B = A$ and for all $a_1, a_2 \in A$, $(a_1, a_2) \in R$ or $(a_2, a_1) \in R$,
            \end{itemize}
        \end{properties}
            
        \begin{definition} 
            An \textit{equivalence relation} on a set $A$ is a reflexive, symmetric and transitive binary relation $\sim \ \subseteq A\times A$. The relation $\sim$ parititons $A$ into equivalence classes $[a]$ where for each $a \in A$, we define $[a] = \{a' \in A: a \sim a'\}$.
        \end{definition}
        
    \section{Functions and Correspondences}
        \begin{definition} 
            Let $X$ and $Y$ be non-empty sets. A \textit{function} from $X$ to $Y$ is a functional left-total binary relation $f \subseteq X\times Y$, which we denote as $f:X\rightarrow Y$. We denote the set of functions from $X$ to $Y$ as $Y^X$.
        \end{definition}
            
        Let $f:X\rightarrow Y$ be a function. Since $f$ is functional, for each $x \in X$ we can denote by $f(x)$ the unique element in $Y$ such that $(x, f(x)) \in f$. The \textit{image} of $f$ is the set $\{f(x): x \in X\}$ which we denote as $\Image(f)$. If $Y = X$ then the function that maps each element of $X$ to itself acts as an identity under composition, we call this the \textit{identity function} and denote it as $\id_X:X\rightarrow X$.
        
        \begin{definition} 
            A \textit{bijection} is an injective and surjective function. Bijections form a group under composition.
        \end{definition}
    
        \begin{definition} 
            A function $f:\mathbb{R}\rightarrow\mathbb{R}$ is \textit{strictly increasing} if for each $x, y \in \mathbb{R}$ we have $x < y$ if and only if $f(x) < f(y)$.
        \end{definition}

        \begin{proposition} 
            Surjective strictly increasing functions form a group under composition.
            
            \begin{proof}
                The identity function $\id_{\mathbb{R}}$ is a surjective increasing function giving the existence of an identity.
            
                Let $f, g:\mathbb{R}\rightarrow\mathbb{R}$ be surjective increasing functions.
                
                If $x, y \in \mathbb{R}$ such that $f(x) = f(y)$ then we must have $x = y$, making $f$ injective. Therefore $f^{-1}$ exists and is also a bijection.
                
                For each $x, y \in \mathbb{R}$ we have $x < y$ if and only if $f(x) < f(y)$, which gives $f^{-1}(x) < f^{-1}(y)$ if and only if $x < y$. Therefore $f^{-1}$ is also a surjective increasing function.
                
                Now $g \circ f$ is surjective and for each $x, y \in \mathbb{R}$ we have $x < y$ if and only if $f(x) < f(y)$ and $f(x) < f(y)$ if and only if $g(f(x)) < g(f(y))$. Therefore $x < y$ if and only if $(g \circ f)(x) < (g \circ f)(y)$ making $g \circ f$ strictly increasing.
            \end{proof}
        \end{proposition}
        
        \begin{definition} 
            A \textit{positive linear transformation} is a function $f:\mathbb{R}\rightarrow\mathbb{R}$ where for some $\alpha \in \mathbb{R}^+$ and $\beta \in \mathbb{R}$, $f(x) = \alpha x + \beta$ for all $x \in \mathbb{R}$.
        \end{definition}
            
        \begin{proposition} 
            Positive linear transformations form a subgroup of the surjective strictly increasing functions under composition.
            
            \begin{proof}
                A positive linear transformation is differentiable and its derivative is positive making it a strictly increasing function.
            
                The identity function which maps real numbers to themselves is a positive linear transformation giving the existence of an identity.
                
                Let $\alpha, \gamma \in \mathbb{R}^+$, $\beta, \delta \in \mathbb{R}$ and $f, g:\mathbb{R}\rightarrow\mathbb{R}$ be positive linear transformations where for each $x \in \mathbb{R}$ we have $f(x) = \alpha x + \beta$ and $g(x) = \gamma x + \delta$.
            
                For each $x \in \mathbb{R}$ we have $f^{-1}(x) = \frac{x - \beta}{\alpha} = \frac{1}{\alpha}x - \frac{\beta}{\alpha}$ and $(g \circ f)(x) = \gamma(\alpha x + \beta) + \delta = \alpha\gamma x + (\beta\gamma + \delta)$ which are both positive linear transformations.
            \end{proof}
        \end{proposition}
    
        \begin{definition} 
            A \textit{correspondence} between two sets $X$ and $Y$ is a left-total binary relation $c \subseteq X\times Y$. For simplicity, we can treat a correspondence $c \subseteq X\times Y$ as a function $c:X\rightarrow \mathcal{P}(Y)^*$ where $\mathcal{P}(Y)^* = \mathcal{P}(Y)-\{\varnothing \}$ and for each $x \in X$ we have $c(x) = \{y \in Y: (x,y) \in c\}$.
        \end{definition}
            
        \begin{definition} 
            A \textit{fixed point} of a correspondence $c:X\rightarrow\mathcal{P}(Y)^*$ is an element $x \in X$ such that $x \in c(x)$.
        \end{definition}
        
        We now state Kakutani's \cite{Kakutani} fixed point theorem without proof which will be referred to when sketching out how to prove the existence of a Nash equilibrium in the Chapter \ref{chap:games}.
        
        \begin{theorem} 
            If $X$ is a non-empty, compact and convex subset of $\mathbb{R}^d$ and $c:X\rightarrow\mathcal{P}(X)^*$ is a correspondence on $X$ with a closed graph and the property that $c(x)$ is convex for all $x \in X$. Then there exists a fixed point $x \in X$ such that $x \in c(x)$.
        \end{theorem}
        
    \section{Preference Relations}
        It is common in economics to require an agent to assign preferences to a set of alternatives $A$ which are captured with the notion of a preference relation.
        
        \begin{definition} 
            A \textit{preference relation} or \textit{total preorder} for an agent over a set $A$ is a reflexive, transitive and total binary relation $\precsim \subseteq A\times A$. If $(a_1, a_2) \in \precsim$ we write $a_1 \precsim a_2$ and say that $a_2$ is \textit{weakly preferred to} or \textit{at least as good as} $a_1$.
        \end{definition}
        
        \begin{definition} 
            Two preference relations $\precsim_1 \subseteq A \times A$ and $\precsim_2 \subseteq B \times B$ are \textit{isomorphic} if there exists a bijection $f:A\rightarrow B$ such that $a_1 \precsim_1 a_2$ if and only if $f(a_1) \precsim_2 f(a_2)$.
        \end{definition}
        
    \section{Relations and Matchings}
	Let $N$ and $\{A_i: i \in N\}$ be non-empty sets. To simplify notation for each $i \in N$ we denote $A_{-i}$ as the set $\times_{j \in N-\{i\}} A_j$, and for $a_i \in A_i$ and $a_{-i} = (a_1, ..., a_{i-1}, a_{i+1}, ...) \in A_{-i}$ we denote $(a_i, a_{-i})$ as $(a_1, ..., a_{i-i}, a_i, a_{i+1}, ...) \in \times_{i \in N} A_i$.
    
    \begin{definition} 
        Let $N$ and $\{A_i: i \in N\}$ be non-empty sets. A \textit{relation} $R$ on  $\{A_i: i \in N\}$ is a subset of $\times_{i \in N} A_i$.
    \end{definition}

        \begin{properties} 
            A relation $R \subseteq \times_{i \in N} A_i$ is; 
            \begin{itemize}
                \item $i-unique$ if for all $a_i \in A_i$ and $a_{-i}, a_{-i}' \in A_{-i}$, $(a_i, a_{-i}), (a_i, a_{-i}') \in R$ implies $a_{-i} = a_{-i}'$.
                \item $i-total$ if for all $a_i \in A_i$ there exists $a_{-i} \in A_{-i}$ such that $(a_i, a_{-i}) \in R$.
            \end{itemize}
        \end{properties}

        \begin{definition} 
            A \textit{matching} $M$ of $n \geq 2$ sets $A_1, ..., A_n$ is a relation which is $i$-total and $i$-unique for all $1 \leq i \leq n$. 
        \end{definition}
        
        A matching of $A_1, ..., A_n$ requires $|A_1| = ... = |A_n|$ and takes the form:
        \[\{(a_1, ..., a_n) : \text{each } a_i \text{ appears in exactly one $n$-tuple}\}.\]
        
        For $1 \leq i, j \leq n$ a matching $M$ induces a unique bijection $M_{ij}:A_i\rightarrow A_j$ which maps elements $a_i \in A_i$ to the unique element $a_j \in A_j$ such that $a_i$ and $a_j$ are elements of an n-tuple in $M$. $M_{ii}$ is the identity bijection for all $1 \leq i \leq n$ and $M_{kl} \circ (M_{jk} \circ M_{ij}) = (M_{kl} \circ M_{jk}) \circ M_{ij}$ for all $1 \leq i, j, k, l \leq n$.
        
    \section{Groupoids}
        \begin{definition} 
            A \textit{groupoid} consists of
            \begin{enumerate}
                \item A non-empty set of objects $C$, and
                \item For each pair of objects $X, Y \in C$ a set of bijective morphisms from $X$ to $Y$ denoted as $C(X, Y)$.
            \end{enumerate}
            
            Satisfying
            \begin{itemize}
                \item For every object $X \in C$ there exists an element $\id_X \in C(X, X)$, 
                \item For every triple of objects $X, Y, Z \in C$ there exists a composite function $\circ:C(Y, Z)\times C(X, Y)\rightarrow C(X, Z)$, and
                \item For every pair of objects $X, Y \in C$ there exists an inverse function $\inv:C(X, Y)\rightarrow C(Y, X)$.
            \end{itemize}
            
            Also satsifying the additional properties that for each $W, X, Y, X \in C$ and $f \in \bij(W, X), g \in \bij(X, Y)$, and $h \in \bij(Y, Z)$
            \begin{itemize}
                \item $f \circ \id_W = f = \id_X \circ f$,
                \item $h \circ (g \circ f) = (h \circ g) \circ f$, and
                \item $f \circ \inv f = \id_Y$ and $\inv f \circ f = \id_X$.
            \end{itemize}
        \end{definition}

    \chapter{Normal-Form Games} \label{chap:games}
    A normal-form game represents a situation where at least two players simultaneously select a strategy from a set of possible strategies. A combination of players' strategies is called a pure strategy profile, and each player has a utility or payoff value specified for every possible pure strategy profile. 
    
    We assume that all games have common knowledge, which is to say that each player knows not only their own strategies and payoffs, but also those of their opponents, knows that each of their opponents know this, ad infinitum. Each player aims to pick their strategy, possibly randomly, to maximise their utility, knowing that each of their opponents are doing the same.
    
    In this chapter we formally define what a normal-form game is and give a brief overview of several concepts which will play a central role in our analysis for the following chapters. For a more thorough coverage of this material see Myerson \cite{Myerson}.
    
    \section{Definition of a Game}
        \begin{definition} 
            A \textit{normal-form game} consists of a set $N$ of at least two players, where each player $i \in N$ has a non-empty set $A_i$ of strategies and a utility or payoff function $u_i:A\rightarrow\mathbb{R}$ where $A = \times_{i \in N}A_i$ is the \textit{pure strategy space} of the game. Elements of $A$ are typically referred to as \textit{pure strategy profiles}. We denote such a game as the triple $\Gamma = (N, A, u)$.
        \end{definition}
            
        A normal-form game is finite if the player and strategy sets are finite, since we only concern ourselves with finite normal-form games we simply refer to them as games. For each $i \in N$ we denote their number of strategies $|A_i|$ as $d_i$.
        
        \begin{proposition}
            For each player $i \in N$, their utility function induces a preference relation $\precsim_i$ over the pure strategy space where $s \precsim_i s'$ if and only if $u_i(s) \leq u_i(s')$ for all $s, s' \in A$.
            
            \begin{proof}
                This follows from $\leq$ being reflexive, transitive and total on $\mathbb{R}$.
            \end{proof}
        \end{proposition}
        
        Although players must select a single strategy when playing a game, they are able to select their strategy using a random process, we call such a strategy a mixed strategy.
        
        \begin{definition} 
            The \textit{mixed strategy space} $\Delta(A_i)$ for player $i$ is the set $\{\sigma_i \in [0,1]^{A_i}: \sum_{s_i \in A_i} \sigma_i(s_i) = 1\}$ of probability distributions over $A_i$. The \textit{mixed strategy space} $\Delta(A)$ for $\Gamma$ is the Cartesian product of the players' mixed strategy spaces $\times_{i \in N} \Delta(A_i)$. Elements of $\Delta(A)$ are typically referred to as \textit{mixed strategy profiles}. For a mixed strategy profile $\sigma = (\sigma_i)_{i \in N} \in \Delta(A)$ and pure strategy $s = (s_i)_{i \in N} \in A$ we denote $\sigma(s)$ as the product $\prod_{i \in N}\sigma_i(s_i)$. 
        \end{definition}
        
        \begin{proposition}
            The mixed strategy space is a convex subset of $\mathbb{R}^{d_1...d_n}$.
            
            \begin{proof}
                Let $\sigma, \sigma' \in \Delta(A)$ and $p \in [0, 1]$. Then 
                \begin{itemize}
                    \item $(p.\sigma_i + (1-p).\sigma_i')(s_i) = p.\sigma_i(s_i) + (1-p).\sigma_i'(s_i) \in [0, 1]$ for all $i \in N$ and $s_i \in A_i$, and
                    \item $\sum_{s_i \in A_i}(p.\sigma_i + (1-p).\sigma_i')(s_i) = p.\sum_{s_i \in A_i}\sigma_i(s_i) + (1-p).\sum_{s_i \in A_i}\sigma_i'(s_i) = 1$ for all $i \in N$.
                \end{itemize}
                
                Therefore $p.\sigma + (1-p).\sigma' \in \Delta(A)$.
            \end{proof}
        \end{proposition}
        
        While we cannot give a definite utility value to a player for a given mixed strategy profile, we can introduce a notion of expected utility by linearly extending the domain of each player's utility function to the mixed strategy space as follows.
        
        \begin{definition} 
            The \textit{expected utility function} for player $i \in N$ is the function $\tilde{u}_i:\Delta(A)\rightarrow\mathbb{R}$ where $\tilde{u}_i(\sigma) = \sum_{s \in A}\sigma(s)u_i(s)$ for all $\sigma \in \Delta(A)$.
        \end{definition}
        
        \begin{proposition}
            For each player $i \in N$, their expected utility function induces a preference relation $\precsim_i$ over the mixed strategy space where $\sigma \precsim_i \sigma'$ if and only if $\tilde{u}_i(\sigma) \leq \tilde{u}_i(\sigma')$ for all $\sigma, \sigma' \in \Delta(A)$.
            
            \begin{proof}
                This follows from $\leq$ being reflexive, transitive and total on $\mathbb{R}$.
            \end{proof}
        \end{proposition}
        
        \begin{definition}
            Let $\sigma, \sigma', \sigma'' \in \Delta(A)$. For each player $i \in N$, their expected utility function is;
            \begin{itemize}
                \item \textit{Continuous} if $\tilde{u}_i(\sigma) \leq \tilde{u}_i(\sigma') \leq \tilde{u}_i(\sigma'')$ implies the existance of $p \in [0,1]$ such that $\tilde{u}_i(p.\sigma + (1-p)\sigma'') = \tilde{u}_i(\sigma')$, and
                \item \textit{Independent} if $\tilde{u}_i(\sigma) \leq \tilde{u}_i(\sigma')$ implies $\tilde{u}_i(p.\sigma + (1-p).\sigma'') \leq \tilde{u}_i(p.\sigma' + (1-p).\sigma'')$ for all $p \in [0,1]$.
            \end{itemize}
        \end{definition}
        
        \begin{proposition}
            For each player $i \in N$, their utility function is continuous and independent.
            
            \begin{proof}
                Suppose $\tilde{u}_i(\sigma) \leq \tilde{u}_i(\sigma') \leq \tilde{u}_i(\sigma'')$. If $\tilde{u}_i(\sigma) = \tilde{u}_i(\sigma'')$ then our continuity condition is satisfied for all $p \in [0,1]$. Now suppose $\tilde{u}_i(\sigma) \neq \tilde{u}_i(\sigma'')$. Then for continuity we need $p \in [0, 1]$ such that,
                \[ \tilde{u}_i(p.\sigma + (1-p).\sigma'') = p.\tilde{u}_i(\sigma) + (1-p).\tilde{u}_i(\sigma'') = \tilde{u}_i(\sigma) \]
                
                Since this is satisfied for $p = \frac{\tilde{u}_i(\sigma') - \tilde{u}_i(\sigma'')}{\tilde{u}_i(\sigma) - \tilde{u}_i(\sigma'')} \in [0,1]$, $\tilde{u}_i$ is continuous.
                
                Now suppose $\tilde{u}_i(\sigma) \leq \tilde{u}_i(\sigma')$. Then for each $p \in [0, 1]$ we have $\tilde{u}_i(p.\sigma + (1-p).\sigma'') = p\tilde{u}_i(\sigma) + (1-p).\tilde{u}_i(\sigma'') \leq p\tilde{u}_i(\sigma') + (1-p).\tilde{u}_i(\sigma'') = \tilde{u}_i(p.\sigma' + (1-p).\sigma'')$.
                
                Therefore $\tilde{u}_i$ is independent.
            \end{proof}
        \end{proposition}
        
        Totality, transitivity, continuity and independence are the assumptions von Neumann and Morgenstern \cite{VNM} require for players' preferences over simple lotteries in order for such an expected utility function to exist in their famous expected utility theorem.
     
    \section{Table Representation of a Game}
        It can be convenient to represent each player's utility values in a single $n$-dimensional table, where each dimension corresponds to the strategy choices of one of the players, and each element of of the table corresponds to a unique pure strategy profile, which contains the utility values for each of the players. To make this clearer we give a couple of examples below.
        
        \begin{example}
            Below is the famous Prisoner's Dilemma game, the row and columns correspond to the strategies of players 1 and 2 respectively.
            \begin{center}
                \begin{game}{2}{2}
                    \> $d$   \> $c$\\
                $d$ \> $2,2$ \> $1,4$\\
                $c$ \> $4,1$ \> $3,3$
                \end{game} 
            \end{center}
            
            \[ A_1 = A_2 = \{d, c\}, A = \{(d,d), (d,c), (c,d), (c,c)\} \]
            Let $s = (d, c) \in A$ and $\sigma = ((0.2, 0.8), (0.5, 0.5)) \in \Delta(A)$, then $u_2(s) = 4$ and 
            \[ \tilde{u}_2(\sigma) = (0.2)(0.5)2 + (0.2)(0.5)4 + (0.8)(0.5)1 + (0.8)(0.5)3 = 2.2 \]
        \end{example}

        \begin{example}
            Below is a $2\times 2 \times 2$ game where each matrix corresponds to one of player 3's strategies, and the row and columns correspond to the strategies of players 1 and 2 respectively.
            \begin{center}
                \begin{game}{2}{2}[$c_1$]
                            \> $b_1$    \> $b_2$ \\
                    $a_1$   \> $1,1,1$  \> $2,3,2$ \\
                    $a_2$   \> $3,2,2$  \> $4,4,5$
                \end{game}
                \hspace*{10mm} 
                \begin{game}{2}{2}[$c_2$]
                            \> $b_1$     \> $b_2$ \\
                    $a_1$   \> $2,2,3$ \> $5,4,4$ \\
                    $a_2$   \> $4,5,4$ \> $6,6,6$
                \end{game}
            \end{center}
            
            \[ A_1 = \{a_1, a_2\}, A_2 = \{b_1, b_2\}, A_3 = \{c_1, c_2\} \]
            Let $s = (a_2, b_1, c_2) \in A$, then $u_2(s) = 5$.
        \end{example}
        
    \section{Dominant Strategies}
        While there is no obvious way for players to give a preference ordering of their strategy set, we can introduce a weaker ordering notion with the idea of a player's strategy dominating another one. 
        
        \begin{definition} 
            Let $s_i, s_i' \in A_i$. $s_i$ \textit{strictly dominates} $s_i'$ if $u_i(s_i, s_{-i}) > u_i(s_i', s_{-i})$ for all $s_{-i} \in A_{-i}$, which we denote as $s_i \succ_i s_i'$.
        \end{definition}
        
        \begin{definition} 
            Let $s_i \in A_i$ and $\sigma_i \in \Delta(A_i)$. $\sigma_i$ \textit{strictly dominates} $s_i$ if $u_i(\sigma_i, \sigma_{-i}) > u_i(s_i, \sigma_{-i})$ for all $\sigma_{-i} \in \Delta(A_{-i})$, which we denote as $\sigma_i \succ_i s_i$.
        \end{definition}
        
    \section{Best Response Correspondences}      
        A strategy $s_i \in A_i$ for player $i \in N$ is a \textit{pure best response} to $s_{-i} \in A_{-i}$ if $u_i(s_i, s_{-i}) \geq u_i(s_i', s_{-i})$ for all $s_i' \in A_i$. For each player $i$ and each strategy $s_{-i} \in A_{-i}$, player $i \in N$ has a non-empty set of best responses. 
    
        \begin{definition} 
            Player $i$'s \textit{pure best response correspondence} is the correspondence $b_i:A_{-i}\rightarrow \mathcal{P}(A_i)^*$ where $b_i(s_{-i}) = \argmax_{s_i \in A_i} u_i(s_i, s_{-i})$ for all $s_{-i} \in A_{-i}$. The game's \textit{pure best response correspondence} is the correspondence $b_{\Gamma}:A\rightarrow\mathcal{P}(A)^*$ where $b(s) = \{s' \in A: s'_i \in b_i(s_{-i}) \forall i \in N\}$ for all $s \in A$.
        \end{definition}
        
        A mixed strategy $\sigma_i \in \Delta(A_i)$ for player $i \in N$ is a best response to $\sigma_{-i} \in \Delta(A_{-i})$ if $\tilde{u}_i(\sigma_i, \sigma_{-i}) \geq \tilde{u}_i(\sigma_i', \sigma_{-i})$ for all $\sigma_i' \in \Delta(A_i)$. For each player $i$ and each strategy $\sigma_{-i} \in \Delta(A_{-i})$, player $i \in N$ has a non-empty set of best responses. 
                    
        \begin{definition} 
            Player $i$'s \textit{best response correspondence} is the correspondence $\tilde{b}_i:\Delta(A_{-i})\rightarrow\mathcal{P}(\Delta(A_i))^*$ where for each $\sigma_{-i} \in \Delta(A_{-i})$, $\tilde{b}_i(\sigma_{-i}) = \argmax_{\sigma_i \in \Delta(A_i)} \tilde{u}_i(\sigma_i, \sigma_{-i})$. The game's \textit{best response correspondence} is the correspondence $\tilde{b}:\Delta(A)\rightarrow\mathcal{P}(\Delta(A))^*$ where for each $\sigma \in \Delta(A)$, $\tilde{b}(\sigma) = \{\sigma' \in \Delta(A): \sigma'_i \in \tilde{b}_i(\sigma_{-i}) \ \forall \ i \in N\}$.
        \end{definition}
        
    \section{Nash Equilibria}
        The most central equilibrium concept in game theory is that of Nash equilibria. They represent strategy profiles where no player can deviate from their specified strategy to increase their payoff, given the strategies being played by their opponents.
        
        \begin{definition}
            A strategy profile $s \in A$ is a \textit{pure strategy Nash equilibrium} if for each player $i$, there does not exist an alternative strategy $s_i' \in A_i$ such that $u_i(s_i', s_{-i}) > u_i(s_i, s_{-i})$. 
        \end{definition}
        
        \begin{definition}
            A strategy profile $\sigma \in \Delta(A)$ is a \textit{Nash equilibrium} if for each player $i \in N$, there does not exist an alternative strategy $\sigma_i' \in \Delta(A_i)$ such that $\tilde{u}_i(\sigma_i', \sigma_{-i}) > \tilde{u}_i(\sigma_i, \sigma_{-i})$. 
        \end{definition}
        
        It is fairly straightforward to see that a mixed strategy profile $\sigma \in \Delta(A)$ is a Nash equilibrium if and only if it is a fixed point of the best response correspondence $\tilde{b}$.
        
        We now state Nash's \cite{NashNCG} famous theorem without proof which gives the eixstence of a Nash equilibria. 
        
        \begin{theorem} 
            Every game has a Nash equilibrium.
        \end{theorem}
        
        One way to prove this is to show that the mixed strategy space of the game is a non-empty, compact and convex subset of $\mathbb{R}^{d_1...d_n}$, that $\tilde{b}$ has a closed graph and that $\tilde{b}(\sigma)$ is convex for all $\sigma \in \Delta(A)$. Existence then follows from Kakutani's fixed point theorem. 

    \chapter{Equivalence of Games} \label{chap:equivalence}
    It is often useful to know whether or not two games are strategically equivalent, for example to determine when two seemingly different situations have the same strategic structure or if a game is symmetric under any notion of the word. 
    
    \section{Game Bijections}
        We use the definition of a game bijection given by Gabarró \textit{et al.} \cite{IsoComplexity}.
    
        \begin{definition} 
            Let $\Gamma_1 = (N, A, u)$ and $\Gamma_2 = (M, B, v)$ be two games. A \textit{bijection} between $\Gamma_1$ and $\Gamma_2$ consists of a bijection $\pi:N\rightarrow M$ and for each player $i \in N$, a bijection $\tau_i:A_i\rightarrow B_{\pi(i)}$. We denote the set of all bijections between $\Gamma_1$ and $\Gamma_2$ as $\bij(\Gamma_1, \Gamma_2)$.
        \end{definition}
            
        For $g = (\pi; (\tau_i)_{i \in N}) \in \bij(\Gamma_1, \Gamma_2)$, $i \in N$, $s_i \in A_i$ and $s \in A$ we use similar notation to that introduced by Stein \textit{et al.} \cite{NoahXE} for symmetric games by denoting $g.i$ as $\pi(i) \in M$, $g.s_i$ as $\tau_i(s_i) \in B_{\pi(i)}$, and $g.s$ as $(\tau_{\pi^{-1}(j)}(s_{\pi^{-1}(j)}))_{j \in M} \in B$ giving $(g.s)_{g.i} = \tau_i(s_i)$.
            
        \begin{definition} 
            For $g = (\pi; (\tau_i)_{i \in N}) \in \bij(\Gamma_1, \Gamma_2)$ and $h = (\eta; (\phi_j)_{j \in M}) \in \bij(\Gamma_2, \Gamma_3)$, their \textit{composite}, denoted $h\circ g$, is $(\eta\circ\pi; (\phi_{\pi(i)}\circ\tau_i)_{i \in N}) \in \bij(\Gamma_1, \Gamma_3)$ giving $(h\circ g).s = h.(g.s)$ for each $s \in A$, and the \textit{inverse} of $g$, denoted $g^{-1}$, is $(\pi^{-1}; (\tau^{-1}_{\pi^{-1}(j)})_{j \in M}) \in \bij(\Gamma_2, \Gamma_1)$.
        \end{definition}
            
        Furthermore, each $g \in \bij(\Gamma_1, \Gamma_2)$ induces $v_{g.i}\circ g \in \mathbb{R}^{A_i}$ where $(v_{g.i}\circ g)(s) = v_{g.i}(g.s)$ for all $s \in A$, and each $g \in \bij(\Gamma_1, \Gamma_2)$ and $\sigma \in \Delta(A)$ induces $\sigma\circ g^{-1} \in \Delta(B)$ where $(\sigma\circ g^{-1})(t) = \sigma(g^{-1}.t)$ for each $t \in B$, giving us $(\sigma\circ g^{-1})(g.s) = \sigma(s)$. When convenient we denote $\sigma \circ g^{-1}$ as $g.\sigma$.
        
        \begin{proposition}
            Games and game bijections are a groupoid.
            
            \begin{proof}
                Let $\Gamma_1 = (N, A, u)$, $\Gamma_2 = (M, B, v)$ and $f = (\pi; (\tau_i)_{i \in N}) \in \bij(\Gamma_1, \Gamma_2)$. Then,
                \begin{align*}
                    f \circ \text{id}_{\Gamma_1} &= (\pi; (\tau_i)_{i \in N}) \circ (\text{id}_N; (\text{id}_{A_i})_{i \in N}) \\
                          &= (\pi \circ \text{id}_N; (\tau_i \circ \text{id}_{A_i})_{i \in N}) \\
                          &= (\pi; (\tau_i)_{i \in N}) = f, 
                \end{align*}
                \begin{align*}
                    \text{id}_{\Gamma_2} \circ f &= (\text{id}_M; (\text{id}_{B_j})_{j \in M}) \circ (\pi; (\tau_i)_{i \in N}) \\
                          &= (\text{id}_M \circ \pi; (\text{id}_{B_{\pi(i)}} \circ \tau_i)_{i \in N}) \\
                          &= (\pi; (\tau_i)_{i \in N}) = f, 
                \end{align*}
                \begin{align*}
                    f \circ f^{-1} &= (\pi; (\tau_i)_{i \in N}) \circ (\pi^{-1}; (\tau^{-1}_{\pi^{-1}(j)})_{j \in M}) \\
                        &= (\pi \circ \pi^{-1}; (\tau_{\pi^{-1}(j)} \circ \tau^{-1}_{\pi^{-1}(j)})_{j \in M}) \\
                        &= (\text{id}_M, (\text{id}_{B_j})_{j \in M}) = \text{id}_{\Gamma_2}, \text{ and} 
                \end{align*}
                \begin{align*}
                    f^{-1} \circ f &= (\pi^{-1}; (\tau^{-1}_{\pi^{-1}(j)})_{j \in M}) \circ (\pi; (\tau_i)_{i \in N}) \\
                        &= (\pi^{-1} \circ \pi; (\tau^{-1}_{\pi^{-1}(\pi(i))} \circ \tau_i)_{i \in N}) \\
                        &= (\text{id}_N, (\text{id}_{A_i})_{i \in N}) = \text{id}_{\Gamma_1}.
                \end{align*}
                
                Now let $\Gamma_3 = (P, C, w)$, $\Gamma_4 = (Q, D, x)$,  $g = (\eta; (\phi_j)_{j \in M}) \in \bij(\Gamma_2, \Gamma_3)$, and $h = (\xi ; (\lambda_k)_{k \in P}) \in \bij(\Gamma_3, \Gamma_4)$. Then, 
                \begin{align*}
                    h \circ (g \circ f) &= (\xi ; (\lambda_k)_{k \in P}) \circ ( (\eta; (\phi_j)_{j \in M}) \circ (\pi; (\tau_i)_{i \in N}) )\\
                        &= (\xi ; (\lambda_k)_{k \in P}) \circ (\eta \circ \pi; (\phi_{\pi(i)} \circ \tau_i)_{i \in N}) \\
                        &= (\xi \circ (\eta \circ \pi); (\lambda_{(\eta \circ \pi)(i)} \circ (\phi_{\pi(i)} \circ \tau_i))_{i \in N}) \\
                        &= (\xi \circ \eta \circ \pi; (\lambda_{(\eta \circ \pi)(i)} \circ \phi_{\pi(i)} \circ \tau_i)_{i \in N}) \text{, and}
                \end{align*}
                \begin{align*}
                    (h \circ g) \circ f &= ( (\xi ; (\lambda_k)_{k \in P}) \circ (\eta; (\phi_j)_{j \in M}) ) \circ (\pi; (\tau_i)_{i \in N}) \\
                        &= (\xi \circ \eta; (\lambda_{\eta(j)} \circ \phi_j)_{j \in M}) \circ (\pi; (\tau_i)_{i \in N}) \\
                        &= ((\xi \circ \eta) \circ \pi; ((\lambda_{\eta(\pi(i))} \circ \phi_{\pi(i)}) \circ \tau_i)_{i \in N}) \\
                        &= (\xi \circ \eta \circ \pi; (\lambda_{(\eta \circ \pi)(i)} \circ \phi_{\pi(i)} \circ \tau_i)_{i \in N}).
                \end{align*}
         
            \end{proof}
        \end{proposition}

    \section {Game Isomorphisms}
        To distinguish whether or not two games have the same structure we need the concept of structure preserving game bijections, called game isomorphisms. A number of possible isomorphisms are possible, depending on how much structure we want preserved.
        
        We begin with the simplest notion of a game isomorphism which was introduced by Nash \cite{NashNCG}, it requires the preservation of players' utility values.
        
        \begin{definition} 
             A bijection $g \in \bij(\Gamma_1, \Gamma_2)$ is a \textit{strict game isomorphism} if $u_i = v_{g.i}\circ g$ for all $i \in N$. We denote the set of strict game isomorphisms between $\Gamma_1$ and $\Gamma_2$ as $\isom(\Gamma_1, \Gamma_2)$.
        \end{definition}

        \begin{example}
            Two games with the same strategic structure.
            
            \begin{center}
                \begin{game}{2}{2}
                            \> $b_1$  \> $b_2$ \\
                    $a_1$   \> $1,8$  \> $2,7$ \\
                    $a_2$   \> $3,6$  \> $4,5$
                \end{game}
                \hspace*{10mm} 
                \begin{game}{2}{2}
                            \> $d_1$ \> $d_2$ \\
                    $c_1$   \> $6,3$ \> $8,1$ \\
                    $c_2$   \> $5,4$ \> $7,2$ \\
                \end{game} 
            \end{center}
        
            Let $\pi = \bigl(\begin{smallmatrix} 1 & 2 \\ 2 & 1 \end{smallmatrix}\bigr), \tau_1 = \bigl(\begin{smallmatrix} a_1 & a_2 \\ d_2 & d_1 \end{smallmatrix}\bigr)$, and $\tau_2 = \bigl(\begin{smallmatrix} b_1 & b_2 \\ c_1 & c_2 \end{smallmatrix}\bigr)$, then $g = (\pi; \tau_1, \tau_2)$ is a strict game isomorphism. For example 
            \[ u_1(a_1, b_2) = v_{\pi(1)}(\tau_2(b_2), \tau_1(a_1)) = v_2(c_2, d_2). \]
        \end{example}
        
        Our second type of game isomorphism (see Gabarro \textit{et al.} \cite{IsoComplexity}) requires the preservation of players' preferences over the the pure strategy spaces.
        
        \begin{definition} 
            A bijection $g \in \bij(\Gamma_1, \Gamma_2)$ is an \textit{ordinal game isomorphism} if for each $i \in N$ we have $u_i(s) \leq u_i(s')$ if and only if $v_{g.i}(g.s) \leq v_{g.i}(g.s')$ for all $s, s' \in A$. We denote the set of ordinal game isomorphisms between $\Gamma_1$ and $\Gamma_2$ as $\isom_o(\Gamma_1, \Gamma_2)$.
        \end{definition}
    
        Our third type of game isomorphism (see Harsanyi and Selten \cite{HarsanyiSelten}) requires the preservation of players' preferences over the mixed strategy spaces.
        
        \begin{definition} 
            A bijection $g \in \bij(\Gamma_1, \Gamma_2)$ is a \textit{cardinal game isomorphism} if for each $i \in N$ we have $\tilde{u}_i(\sigma) \leq \tilde{u}_i(\sigma')$ if and only if $\tilde{v}_{g.i}(g.\sigma) \leq \tilde{v}_{g.i}(g.\sigma')$ for all $\sigma, \sigma' \in \Delta(A)$. We denote the set of cardinal game isomorphisms as $\isom_c(\Gamma_1, \Gamma_2)$.
        \end{definition}
            
    \section{Properties of Game Isomorphisms} 
        \begin{proposition} 
            Games and strict, ordinal or cardinal game isomorphisms are a groupoid.
            
            \begin{proof}
                We just prove this for ordinal game isomorphisms as the other proofs are similar.
                
                If $g \in \isom_o(\Gamma_1, \Gamma_2)$ then we have $u_i(s) \leq u_i(s')$ if and only if $v_{g.i}(g.s) \leq v_{g.i}(g.s')$, which is equivalent to $u_{g^{-1}.j}(g^{-1}.t) \leq u_{g^{-1}.j}(g^{-1}.t')$ if and only if $v_j(t) \leq v_j(t')$. Therefore $g^{-1} \in \isom_o(\Gamma_2, \Gamma_1)$.
                
                Let $\Gamma_1 = (N, A, u)$, $\Gamma_2 = (M, B, v)$ and $\Gamma_3 = (P, C, w)$ such that there exists $g \in \isom_o(\Gamma_1, \Gamma_2)$ and $h \in \isom_o(\Gamma_2, \Gamma_3)$. 
                
                Then for each $i \in N$ and $s, s' \in A$ we have, $u_i(s) \leq u_i(s')$ if and only if $v_{g.i}(g.s) \leq v_{g.i}(g.s')$, and $v_{g.i}(g.s) \leq v_{g.i}(g.s')$ if and only if $w_{h.(g.i)}(h.(g.s)) \leq w_{h.(g.i)}(h.(g.s'))$.
                
                Therefore $u_i(s) \leq u_i(s')$ if and only if $w_{(h\circ g).i}((h\circ g).s) \leq w_{(h\circ g).i}((h\circ g).s')$, giving $h \circ g \in isom_o(\Gamma_1, \Gamma_3)$.
                
                Our other conditions follow from games and game bijections being a groupoid.
            \end{proof}
        \end{proposition}
            
        \begin{proposition} \label{ordisoprop}
            $g \in \bij(\Gamma_1, \Gamma_2)$ is an ordinal game isomorphism if and only if for each player $i \in N$ there exists a surjective strictly increasing function $\alpha_i \in \mathbb{R}^{\mathbb{R}}$ such that $\alpha_i \circ u_i = v_{g.i}\circ g$.
            
            \begin{proof}
                Suppose there exists $g \in \bij(\Gamma_1, \Gamma_2)$ and a surjective strictly increasing function $\alpha_i \in \mathbb{R}^{\mathbb{R}}$ such that $\alpha_i \circ u_i = v_{g.i}\circ g$. Then for each $i \in N$ and $s, s' \in A$ such that $u_i(s) \leq u_i(s')$ we have,
                \[ (\alpha \circ u_i)(s) = v_{g.i}(g.s) \leq v_{g.i}(g.s') = (\alpha \circ u_i)(s'). \]
                
                Similarly, for each $j \in M$ and $t, t' \in B$ such that $v_j(t) \leq v_j(t')$, we have,
                \[ (\alpha^{-1} \circ v_j)(t) = u_{g^{-1}.j}(g^{-1}.t) \leq u_{g^{-1}.j}(g^{-1}.t') = (\alpha^{-1} \circ v_j)(t'). \]
                
                Conversely, suppose there exists $g \in \isom_o(\Gamma_1, \Gamma_2)$. 
                
                Let $i \in N$ and $\bar{A}_0 = \argmin_{s \in A} u_i(s)$. 
                
                For each $l \in \mathbb{Z}^+$ let $\bar{A}_l = \argmin_{s \in A - \cup^{l-1}_{j=1}\bar{A}_j} u_i(s)$.
                
                Let $k \in \mathbb{N}$ such that $\bar{A}_k$ is non-empty and $\bar{A}_{k+1}$ is empty.
                
                Let $s_0 \in \bar{A}_0$, ..., and $s_k \in \bar{A}_k$. 
                
                Let $\mu_0 = u_i(s_0)$, ..., $\mu_k = u_i(s_k)$, $\nu_0 = v_{g.i}(g.s_0)$, ..., and $\nu_k = v_{g.i}(g.s_k)$.
                
                Define $\alpha_i \in \mathbb{R}^{\mathbb{R}}$ by 
                \[
                    \alpha_i(x) =
                        \begin{cases}
                            \frac{\nu_0}{\mu_0}x  & \text{if } -\infty < x \leq \mu_0 \\
                            \frac{\nu_1 - \nu_0}{\mu_1 - \mu_0}x + \frac{\mu_1\nu_0 - \nu_1\mu_0}{\mu_1 - \mu_0} & \text{if } \mu_0 < x \leq \mu_1 \\
                                 &  \vdots \\
                            \frac{\nu_k - \nu_{k-1}}{\mu_k - \mu_{k-1}}x + \frac{\mu_k\nu_{k-1} - \nu_k\mu_{k-1}}{\mu_k - \mu_{k-1}}  & \text{if } \mu_{k-1} < x \leq \mu_k \\
                            \frac{\nu_k}{\mu_k}x  & \text{if } \mu_k < x < \infty
                        \end{cases}
                \]
                
                Then $\alpha$ is a surjective strictly increasing function with $\alpha_i \circ u_i = v_{g.i}\circ g$.
            \end{proof}
        \end{proposition}
            
        This shows that the players' preferences over the pure strategy spaces are preserved uniquely up to surjective strictly increasing functions. 
        
        We now prove some of the additional properties relating to the pure strategy spaces which are preserved.
            
        \begin{proposition} 
            If $g \in \isom_o(\Gamma_1, \Gamma_2)$ then by \ref{ordisoprop} for each $i \in N$ there exists a surjective strictly increasing function $\alpha_i \in \mathbb{R}^{\mathbb{R}}$ such that $\alpha_i \circ u_i = v_{g.i} \circ g$, and $g$ preserves; 
            
            \begin{enumerate}
                \item The player's pure best response correspondences.
                \item The game's pure best response correspondence.
                \item The pure strategy Nash equilibria.
                \item Dominance of pure strategies over pure strategies.
            \end{enumerate}
            
            \begin{proof}
                \begin{enumerate}
                    \item Let $i \in N$, $s_{-i} \in A_{-i}$ and $s_i \in b_i(s_{-i})$ so that $u_i(s_i, s_{-i}) \leq u_i(s_i', s_{-i})$ for all $s_i' \in A_i$. Then for each $s_i' \in A_i$ we have,
                          \[(\alpha_i \circ u_i)(s_i, s_{-i}) = v_{g.i}(g.(s_i, s_{-i})) \leq v_{g.i}(g.(s_i', s_{-i})) = (\alpha_i \circ u_i)(s_i', s_{-i}).\] 
                          Therefore $g.s_i \in b_{g.i}(g.s_{-i})$.
                        
                          The reverse argument is the same.
                    \item Let $s' \in A$ and $s \in b(s')$ so that $s_i \in b_i(s_{-i}')$ for all $i \in N$. Then $g.s_i \in b_{g.i}(g.s_{-i}')$ for all $i \in N$. Therefore $g.s \in b(g.s')$.
                            
                          The reverse argument is the same.
                    \item Let $s \in A$ be a pure strategy Nash equilibrium. Then $s \in b(s)$ giving us $g.s \in b(g.s)$, making $g.s$ a pure strategy Nash equilibrium.
                            
                          The reverse argument is the same.
                    \item Let $i \in N$ and $s_i, s_i' \in A_i$ such that $u_i(s_i, s_{-i}) \geq u_i(s_i', s_{-i})$ for all $s_{-i} \in A_{-i}$. Then $v_{g.i}(g.(s_i , s_{-i})) \geq v_{g.i}(g.(s_i', s_{-i}))$ for all $s_{-i} \in A_{-i}$. Therefore $g.s_i$ strictly dominates $g.s_i'$.
                            
                          The reverse argument is the same.
                \end{enumerate}
            \end{proof}
        \end{proposition}
        
        Note that the proofs for Propositions \ref{prop:firstnoahproof} and \ref{cardisoprop} were communicated to the author via email by Noah Stein.
        
        \begin{proposition} \label{prop:firstnoahproof}
            Let $\Gamma_1 = (N, A, u)$, $\Gamma_2 = (M, B, v)$, $i \in N$ and $g \in \bij(\Gamma_1, \Gamma_2)$. If there exists a positive linear transformation $\alpha_i \in \mathbb{R}^{\mathbb{R}}$ such that $\alpha_i \circ u_i = v_{g.i} \circ g$, then $\alpha_i \circ \tilde{u}_i = \tilde{v}_{g.i} \circ g$.
            
            \begin{proof}
                Since $\alpha_i$ is a positive linear transformation, there exists $\beta_i \in \mathbb{R}^+$ and $\gamma_i \in \mathbb{R}$ such that for each $s \in A$ we have $\beta_i.u_i(s) + \gamma_i = v_{g.i}(g.s)$. Therefore for each $\sigma \in \Delta(A)$ we have,
            
                \begin{align*} 
                    (\tilde{v}_{g.i}\circ g)(\sigma) = \tilde{v}_{g.i}(\sigma\circ g^{-1}) &= \sum_{t \in B}\sigma\circ g^{-1}(t).v_{g.i}(t) \\
                        &= \sum_{s \in A}\sigma \circ g^{-1}(g.s).v_{g.i}(g.s)   \\
                        &= \sum_{s \in A}\sigma(s).v_{g.i}(g.s)  \\  
                        &= \sum_{s \in A}\sigma(s).(\beta_i.u_i(s) + \gamma_i) \\  
                        &= \beta_i\sum_{s \in A}(\sigma(s).u_i(s)) + \gamma_i\sum_{s \in A}\sigma(s) \\
                        &= \beta_i.\tilde{u}_i(\sigma) + \gamma_i = (\alpha_i \circ \tilde{u}_i)(\sigma)
                \end{align*}
            \end{proof}
        \end{proposition}
    
        \begin{proposition} \label{cardisoprop}
            $g \in \bij(\Gamma_1, \Gamma_2)$ is a cardinal game isomorphism if and only if for each player $i \in N$ there exists a positive linear transformation $\alpha_i:\mathbb{R}\rightarrow\mathbb{R}$ such that $\alpha_i \circ u_i = v_{g.i}\circ g$.
            
            \begin{proof}
                Suppose for each player $i \in N$ there exists a positive linear transformation $\alpha_i \in \mathbb{R}^{\mathbb{R}}$ such that $\alpha_i \circ u_i = v_{g.i}\circ g$. Then for each $i \in N$ and $\sigma, \sigma' \in \Delta(A)$ such that $\tilde{u}_i(\sigma) \leq \tilde{u}_i(\sigma')$ we have,
                \[ \tilde{v}_{g.i}(\sigma \circ g^{-1}) = (\alpha \circ \tilde{u}_i)(\sigma) \leq (\alpha \circ \tilde{u}_i)(\sigma') = \tilde{v}_{g.i}(\sigma' \circ g^{-1}). \]
                
                Similarly, for each $j \in M$ and $\varsigma, \varsigma' \in B$ such that $\tilde{v}_j(\varsigma) \leq \tilde{v}_j(\varsigma')$, we have,
                \[ (\alpha^{-1} \circ \tilde{v}_j)(\varsigma) = \tilde{u}_{g^{-1}.j}(\varsigma \circ g) \leq \tilde{u}_{g^{-1}.j}(\varsigma' \circ g) = (\alpha^{-1} \circ v_j)(\varsigma'). \]
            
                Conversely, suppose $g \in \isom_c(\Gamma_1, \Gamma_2)$ and let $i \in N$, $\underline{s} = \argmin_{s \in A} u_i(s)$, $\bar{s} = \argmax_{s \in A} u_i(s)$. 
                
                Let $\underline{u} = u_i(\underline{s})$, $\bar{u} = u_i(\bar{s})$, $\underline{v} = v_{g.i}(g.\underline{s})$, and $\bar{v} = v_{g.i}(g.\bar{s})$.
                
                Define $\alpha_i \in \mathbb{R}^{\mathbb{R}}$ by
                \[
                    \alpha_i(x) =
                        \begin{cases}
                            x &\text{if } \underline{u} = \bar{u} = 0 \\
                            \frac{\bar{v}}{\bar{u}}x & \text{if } \underline{u} = \bar{u} \neq 0 \\
                            \frac{\bar{v} - \underline{v}}{\bar{u} - \underline{u}}x + \frac{\bar{u}\underline{v} - \bar{v}\underline{u}}{\bar{u} - \underline{u}} & \text {if } \underline{u} \neq \bar{u} \\
                        \end{cases}
                    \text{ for all } x \in \mathbb{R}
                \]
                
                Then $\alpha_i$ is a positive linear transformation with $(\alpha_i \circ u_i)(\underline{s}) = v_{g.i}(g.\underline{s})$ and $(\alpha_i \circ u_i)(\bar{s}) = v_{g.i}(g.\bar{s})$.
                             
                Now let $s \in A$. Since $u_i(\underline{s}) \leq u_i(s) \leq u_i(\bar{s})$ and $\tilde{u}_i$ is continuous, there exists $\sigma \in \Delta(A)$ such that $\sigma(s') = 0$ for all $s' \in A-\{\underline{s}, \bar{s}\}$ and $u_i(s) = \tilde{u}_i(\sigma)$. Therefore,
                \begin{align*}
                    (\alpha_i \circ u_i)(s) &= (\alpha_i \circ \tilde{u}_i)(\sigma) \\
                        &= \sigma(\underline{s})(\alpha_i \circ u_i)(\underline{s}) + \sigma(\bar{s})(\alpha_i \circ u_i)(\bar{s}) \\
                        &= \sigma(\underline{s})v_{g.i}(g.\underline{s}) + \sigma(\bar{s})v_{g.i}(g.\bar{s}) \\
                        &= \tilde{v}_{g.i}(\sigma \circ g^{-1}) \\
                        &= v_{g.i}(g.s)
                \end{align*}
            \end{proof}
        \end{proposition}
        
        This shows that the players' preferences over the mixed strategy spaces are preserved uniquely up to positive linear transformations of their utility functions. 
        
        We now prove some of the additional properties relating to the mixed strategy spaces which are preserved.
            
        \begin{proposition} 
            If $g \in \isom_c(\Gamma_1, \Gamma_2)$ then by \ref{cardisoprop} for each player $i \in N$ there exists $\beta_i \in \mathbb{R}^+$ and $\gamma_i \in \mathbb{R}$ such that $\beta_i.\tilde{u}_i(\sigma) + \gamma_i = \tilde{v}_{g.i}(\sigma \circ g^{-1})$ for all $\sigma \in \Delta(A)$, and $g$ preserves; 
            
            \begin{enumerate}
                \item The player's best response correspondences.
                \item The game's best response correspondence.
                \item The mixed strategy Nash equilibria.
                \item Dominance of mixed strategies over pure strategies.
            \end{enumerate}
            
            \begin{proof}
                \begin{enumerate}
                    \item Let $i \in N$, $\sigma_{-i} \in \Delta(A_{-i})$ and $\sigma_i \in \tilde{b}_i(\sigma_{-i})$ so that $\tilde{u}_i(\sigma_i, \sigma_{-i}) \geq \tilde{u}_i(\sigma_i', \sigma_{-i})$ for all $\sigma_i' \in \Delta(A_i)$.
                            
                          Then $\tilde{v}_{g.i}((\sigma_i, \sigma_{-i}) \circ g^{-1}) \geq \tilde{v}_{g.i}((\sigma_i', \sigma_{-i}) \circ g^{-1})$ for all $\sigma_i' \in \Delta(A_i)$. Therefore $\sigma_i \circ g^{-1} \in \tilde{b}_{g.i}(\sigma_{-i} \circ g^{-1})$.
                            
                          The reverse argument is the same.
                    \item Let $\sigma' \in \Delta(A)$ and $\sigma \in \tilde{b}(\sigma')$, so that $\sigma_i \in \tilde{b}_i(\sigma_{-i}')$ for all $i \in N$. Then $\sigma_i \circ g^{-1} \in \tilde{b}_{g.i}(\sigma_{-i}' \circ g^{-1})$ for all $i \in N$. Therefore $\sigma \circ g^{-1} \in \tilde{b}(\sigma' \circ g^{-1})$.
                            
                          The reverse argument is the same.
                    \item Let $\sigma \in \Delta(A)$ be a Nash equilibrium. Then $\sigma \in \tilde{b}(\sigma)$ giving us $\sigma \circ g^{-1} \in \tilde{b}(\sigma \circ g^{-1})$. Therefore $\sigma \circ g^{-1}$ is a Nash equilibrium.
                            
                          The reverse argument is the same.
                    \item Let $\sigma_i \in \Delta(A_i)$ and $s_i \in A_i$ such that $\tilde{u}_i(\sigma_i, \sigma_{-i}) \geq \tilde{u}_i(s_i, \sigma_{-i})$ for all $\sigma_{-i} \in \Delta(A_{-i})$.
                        
                          Then $\tilde{v}_{g.i}((\sigma_i , \sigma_{-i}) \circ g^{-1}) \geq \tilde{v}_{g.i}((s_i, \sigma_{-i}) \circ g^{-1})$ for all $\sigma_{-i} \in \Delta(A_{-i})$. Therefore  $\sigma \circ g^{-1}$ strictly dominates $g.s$.
                            
                          The reverse argument is the same.
                \end{enumerate}
            \end{proof}
        \end{proposition}
        
    \section{Notions of Equivalence} 
        
        \begin{definition} 
            Let $\Gamma_1$ and $\Gamma_2$ be two games.
            \begin{itemize}
                \item $\Gamma_1$ and $\Gamma_2$ are \textit{strictly equivalent} if $\isom(\Gamma_1, \Gamma_2)$ is non-empty, which we denote as $\Gamma_1 \cong \Gamma_2$,
                \item $\Gamma_1$ and $\Gamma_2$ are \textit{ordinally equivalent} if $\isom_o(\Gamma_1, \Gamma_2)$ is non-empty, which we denote as $\Gamma_1 \cong_o \Gamma_2$, and 
                \item $\Gamma_1$ and $\Gamma_2$ are \textit{cardinally equivalent} if $\isom_c(\Gamma_1, \Gamma_2)$ is non-empty, which we denote as $\Gamma_1 \cong_c \Gamma_2$. 
            \end{itemize}
        \end{definition}
    
        \begin{proposition} 
            $\cong$, $\cong_o$ and $\cong_c$ are equivalence relations.
            
            \begin{proof}
                For each game $\Gamma$ we have $\id_{\Gamma} \in \isom(\Gamma, \Gamma)$, $\id_{\Gamma} \in \isom_o(\Gamma, \Gamma)$ and $\id_{\Gamma} \in \isom_c(\Gamma, \Gamma)$ which gives reflexivity of $\cong$, $\cong_o$ and $\cong_c$.
                
                Our symmetric and transitivity conditions follow from games and each type of game isomorphisms being a groupoid.
            \end{proof}
        \end{proposition}
        
        \begin{definition} 
            Let $\Gamma$ be a game.
            \begin{itemize}
                \item The \textit{strict equivalence class} of $\Gamma$ is the set $\{\Gamma': \Gamma \cong \Gamma'\}$ which we denote as $[\Gamma]$.
                \item The \textit{ordinal equivalence class} of $\Gamma$ is the set $\{\Gamma': \Gamma \cong \Gamma'\}$ which we denote as $[\Gamma]_o$.
                \item The \textit{cardinal equivalence class} of $\Gamma$ is the set $\{\Gamma': \Gamma \cong \Gamma'\}$ which we denote as $[\Gamma]_c$.
            \end{itemize}
        \end{definition}
        
        All games in each equivalence class have the same strategic structure under that notion. We note that for each set of $n$ player games where each player has $d_i$ strategies there are an infinite number of strict and cardinal equivalence classes, and a finite number of ordinal equivalence classes.
        
        Goforth and Robinson \cite{GoforthRobinson} counted 144 ordinal equivalence classes for the two player games where each player has two strategies and their utility function induces a total order over the pure strategy spaces. 
    
        We now show that strict equivalence is a weaker notion of equivalence than cardinal equivalence, and that cardinal equivalence is a weaker notion of equivalnece than ordinal equivalence. 
    
        \begin{proposition} 
            $\isom(\Gamma_1, \Gamma_2) \subseteq \isom_c(\Gamma_1, \Gamma_2) \subseteq \isom_o(\Gamma_1, \Gamma_2)$.
            
            \begin{proof}
                $\isom(\Gamma_1, \Gamma_2) \subseteq \isom_c(\Gamma_1, \Gamma_2)$ since the idenity function $\id_{\mathbb{R}}$ is a positive linear transformation and $\isom_c(\Gamma_1, \Gamma_2) \subseteq \isom_o(\Gamma_1, \Gamma_2)$ since every positive linear transformation is a surjective strictly increasing function.
            \end{proof}
        \end{proposition}
        
        Of course by transitivity of $\subseteq$ this implies that strict equivalence is also a weaker notion of equivalence than cardinal equivalence.

    \chapter{Symmetric Games} \label{chap:symmgames}
    Symmetric games capture the structural notion of a game being fair. Not only is this useful for classifying games but this structure can also be exploited for computational and theoretical purposes. 
    
    In the simplest sense, we would like a symmetric game to represent the same situation regardless of the player. Under any such definition it is necessary for all players to have the same number of strategies so we make this assumption for the remainder of the chapter. 

    \section {Automorphism Group}
        \begin{definition} 
            An isomorphism $g \in \isom(\Gamma, \Gamma)$ is an $automorphism$ or $symmetry$ of $\Gamma$. We denote the set of $automorphisms$ as $\Aut(\Gamma)$.
        \end{definition}
        
        \begin{example}
            The Prisoner's Dilemma
            \begin{center}
                \begin{game}{2}{2}
                          \> $d$    \> $c$ \\
                    $d$   \> $3,3$  \> $1,4$ \\
                    $c$   \> $4,1$  \> $2,2$
                \end{game} 
            \end{center}
            
            \[ Aut(\Gamma) = \{ \bigl(e ; \bigl(\begin{smallmatrix} d & c \\ d & c \end{smallmatrix}\bigr), \bigl(\begin{smallmatrix} d & c \\ d & c \end{smallmatrix}\bigr)\bigr),
                                \bigl((12) ; \bigl(\begin{smallmatrix} d & c \\ d & c \end{smallmatrix}\bigr), \bigl(\begin{smallmatrix} d & c \\ d & c \end{smallmatrix}\bigr)\bigr) \}  \]
        \end{example}
            
        \begin{proposition} 
            The automorphisms of a game form a group under composition.
            
            \begin{proof}
                This follows from the fact that games and strict game isomorphisms are a groupoid.
            \end{proof}
        \end{proposition}
        
        \begin{remark}
            Let $M$ be a matching of the strategy sets. For each player permutation $\pi \in S_N$, $M$ induces a game bijection $(\pi; (M_{i\pi(i)})_{i \in N}) \in \bij(\Gamma, \Gamma)$ which we denote as $M_{\pi}$.
        \end{remark}

    \section {Symmetry Groups}    
        Below we define symmetry groups for games which were introduced by Stein \textit{et al.} \cite{NoahXE}. Symmetry groups make the classification of symmetric games in the next section much clearer.
            
        \begin{definition} 
            A \textit{symmetry group} of a game $\Gamma$ is a subgroup $G$ of the automorphism group $\Aut(\Gamma)$. The \textit{stabiliser subgroup} for player $i$ is the subgroup of automorphisms $\{g \in G: g.i = i\}$ that map player $i$ to itself, which we denote as $G_i$.
        \end{definition}
            
        \begin{properties} 
            Let $G$ be a symmetry group of a game $\Gamma = (N, A, u)$. We say that $G$ is;
            \begin{itemize}
                \item \textit{player transitive} if $G$ acts transitively on $N$, that is for each $i, j \in N$ there exists $g \in G$ such that $g.i = j$,
                \item \textit{player n-transitive} if $G$ acts $n$-transitively on $N$, that is for each $\pi \in S_N$ there exists $g \in G$ such that $g.i = \pi(i)$ for all $i \in N$, and
                \item \textit{strategy trivial} if for each $g \in G_i$, $g.s_i = s_i$ for all $s_i \in A_i$.
            \end{itemize}
        \end{properties}

    \section {Notions of Symmetry}
        In this section we introduce various notions of a game being symmetric. Von Neumann and Morgenstern \cite{VNM} were the first to consider symmetric games, which they defined as follows. 

        \begin{definition} 
            A game $\Gamma$ is \textit{VNM symmetric} if $A_i = A_j$ for all $i, j \in N$, and for each permutation $\pi \in S_N$,  $u_{\pi(i)}(s_1, ..., s_n) = u_i(s_{\pi(1)}, ..., s_{\pi(n)})$ for all $i \in N$ and $(s_1, ..., s_n) \in A$.
        \end{definition}
        
        \begin{example}
            Two Player VNM Symmetric Game
            \begin{center}
                \begin{game}{2}{2}
                          \> $a$    \> $b$ \\
                    $a$   \> $1,1$  \> $3,2$ \\
                    $b$   \> $2,3$  \> $4,4$
                \end{game} 
            \end{center}
            
            For $\pi = \bigl(\begin{smallmatrix} 1 & 2 \\ 2 & 1 \end{smallmatrix}\bigr) \in S_2$ we have
            \begin{align*}
                u_1(a, a) &= u_2(a, a) = 1   &   u_1(a, b) &= u_2(b, a) = 3 \\
                u_1(b, a) &= u_2(a, b) = 2   &   u_1(b, b) &= u_2(b, b) = 4
            \end{align*} 
        \end{example}
        
        While VNM symmetric games are obviously fair, they do not fully capture the notion of fairness, nor do they classify the possible notions of structural fairness that may be present. We now define the classifications of symmetric games which we think best capture the possible structural notions of fairness.
        
        Nash \cite{NashNCG} provided the most general definition of symmetric game as follows.
        
        \begin{definition} 
            A game $\Gamma$ is \textit{symmetric} if its automorphism group is player transitive.
        \end{definition}
        
        Not only does this drop the requirement that strategy sets be equal, it also drops the assumption that we are able to match up the strategy sets such that the game be automorphic under the induced game bijections or that the game be automorphic under every permutation of the players.
        
        Our next notion of a symmetric game was introduced by Stein \textit{et al.} \cite{NoahXE} which strengthens the conditions for a symmetric game.
        
        \begin{definition} 
            A game $\Gamma$ is \textit{standard symmetric} if there exists a player transitive and strategy trivial symmetry group.
        \end{definition}
        
        All symmetric games presented so far have been standard symmetric, below we provide an example of a symmetric game which is not standard symmetric.
        
        \begin{example}
            Matching Pennies
            \begin{center}
                \begin{game}{2}{2}
                          \> $H$    \> $T$ \\
                    $H$   \> $1,-1$  \> $-1,1$ \\
                    $T$   \> $-1,1$  \> $1,-1$
                \end{game} 
            \end{center}
            
            \begin{align*}
                 \Aut(\Gamma) = \{ &\bigl(e ; \bigl(\begin{smallmatrix} H & T \\ H & T \end{smallmatrix}\bigr), \bigl(\begin{smallmatrix} H & T \\ H & T \end{smallmatrix}\bigr)\bigr), 
                                \bigl(e ; \bigl(\begin{smallmatrix} H & T \\ T & H \end{smallmatrix}\bigr), \bigl(\begin{smallmatrix} H & T \\ T & H \end{smallmatrix}\bigr)\bigr), \\
                                &\bigl((12) ; \bigl(\begin{smallmatrix} H & T \\ H & T \end{smallmatrix}\bigr), \bigl(\begin{smallmatrix} H & T \\ T & H \end{smallmatrix}\bigr)\bigr),
                                \bigl((12) ; \bigl(\begin{smallmatrix} H & T \\ T & H \end{smallmatrix}\bigr), \bigl(\begin{smallmatrix} H & T \\ H & T \end{smallmatrix}\bigr)\bigr) \}
            \end{align*}
            Since $\Aut(\Gamma)$ is player transitive, is not strategy trivial and contains no proper subgroups, matching pennies is an example of a symmetric game which is not standard symmetric.
        \end{example}
        
        We now make another stricter definition of a symmetric game.
        
        \begin{definition} 
            A game $\Gamma$ is \textit{fully symmetric} if its automorphism group is player $n$-transitive.
        \end{definition}
        
        Finally, we make the following definitions which partition the set of symmetric games up into different classes based on the kind of structure they contain.
        
        \begin{definition} 
            A game $\Gamma$ is;
            \begin{itemize}
                \item \textit{non-fully non-standard symmetric} if it is symmetric but neither fully or standard symmetric,
                \item \textit{fully non-standard symmetric} if it is fully symmetric and not standard symmetric,
                \item \textit{non-fully standard symmetric} if it is standard symmetric and not fully symmetric, and
                \item \textit{fully standard symmetric} if it is fully symmetric and standard symmetric.
            \end{itemize}
        \end{definition}
        
    \section{Properties of Symmetric Games}
        VNM symmetric games do not have the desirable property that if $\Gamma_1 \cong \Gamma_2$ and $\Gamma_1$ is VNM symmetric then $\Gamma_2$ is also VNM symmetric, which we illustrate with the example below.
        
        \begin{example}
            Consider the two following equivalent games, the first being the prisoners dilemma.
            
            \begin{center}
                \begin{game}{2}{2}[$\Gamma_1$]
                            \> $d$  \> $c$ \\
                    $d$   \> $3,3$  \> $1,4$ \\
                    $c$   \> $4,1$  \> $2,2$
                \end{game}
                \hspace*{10mm} $\cong$ \hspace*{10mm}
                \begin{game}{2}{2}[$\Gamma_2$]
                            \> $a$ \> $b$ \\
                    $a$   \> $1,4$ \> $3,3$ \\
                    $b$   \> $2,2$ \> $4,1$ \\
                \end{game} 
            \end{center}
        
            \[ Isom(\Gamma_1, \Gamma_2) = \{ \bigl(e ; \bigl(\begin{smallmatrix} d & c \\ a & b \end{smallmatrix}\bigr), \bigl(\begin{smallmatrix} d & c \\ b & a \end{smallmatrix}\bigr)\bigr),
                                \bigl((12) ; \bigl(\begin{smallmatrix} d & c \\ b & a \end{smallmatrix}\bigr), \bigl(\begin{smallmatrix} d & c \\ a & b \end{smallmatrix}\bigr)\bigr) \}  \]
                                
            Clearly $\Gamma_1$ is VNM symmetric and $\Gamma_1 \cong \Gamma_2$, but $\Gamma_2$ is not VNM symmetric.
        \end{example}

        We now show that standard symmetric games capture the notion that the game is symmetric under some matching of the strategy sets, we note that we use the argument from Stein \textit{et al.} \cite{NoahXE} to show that a player transitive and strategy trivial symmetry group maps between strategy sets in a canonical manner.
        
        \begin{proposition} 
            A game $\Gamma$ is standard symmetric if and only if there exists a matching of the strategy sets such that each game bijection induced from the player permutations in some transitive subgroup of $S_N$ is an automorphism.
            
            \begin{proof}
                Let $\Gamma$ be a standard symmetric game, $G$ any player transitive and strategy trivial subgroup of $Aut(\Gamma)$, and $g, h \in G$ such that $g.i = h.i = j$. 
                
                For each $g, h \in G$ such that $g.i = h.i = j$ we have $h^{-1} \circ g \in G_i$, strategy triviality of $G$ gives $(h^{-1} \circ g).s_i = s_i$ for all $s_i \in A_i$, so we must have $g.s_i = h.s_i$ for all $s_i \in A_i$. Hence $G$ maps $A_i$ to $A_j$ in a canonical way.
                
                This along with player transitivity of $G$ induces a matching of the strategy sets such that each game bijection induced from the player bijections in some transitive subgroup of $S_N$ is an automorphism.
                
                Conversely, let $\Gamma$ be a game, $H$ a transitive subgroup of $S_N$ and $M$ a matching of the strategy sets such that each game bijection induced from the player permutations in $H$ is an automorphism.
                
                Then $\{M_{\pi} \in \Aut(\Gamma): \pi \in H\}$ is a player transitive and strategy trivial subgroup of $\Aut(\Gamma)$.
            \end{proof}
        \end{proposition}
        
        \begin{proposition}
            If $\Gamma$ is standard symmetric so that there exists a matching $M$ of the strategy sets such that each game bijection induced from the player permutations in some transitive subgroup $H$ of $S_N$ is an automorphism, then for each $i, j \in N$ we have $u_i(s) = u_j(s)$ for all $s \in M$.
            
            \begin{proof}
                This follows from $\{M_{\pi} \in \Aut(\Gamma): \pi \in S_N\}$ being a player transitive symmetry group and that $M_{\pi}.s = s$ for all $s \in M$ and $\pi \in S_N$.
            \end{proof}
        \end{proposition}
        
        This shows that in a standard symmetric game all of the players must assign the same utility to each strategy profile in some matching of the pure strategy space.
        
        \begin{proposition}
            If a game has a player $n-transitive$ and strategy trivial symmetry group then it is fully standard symmetric.
            
            \begin{proof}
                This follows directly from the definition of fully and standard symmetric games.
            \end{proof}
        \end{proposition}
    
        It is not entirely clear whether a fully standard symmetric game implies the existence of a player $n$-transitive and strategy trivial symmetry group, although it would be desirable if it did for a couple of reasons discussed below. We have been unable to show this however, so have left it as a conjecture.
        
        \begin{conjecture} 
            A game is fully standard symmetric if and only if it has a player $n$-transitive and strategy trivial symmetry group. (Note: conjecture is false, see the more recent paper `Notions of Symmetry for Finite Strategic-Form Games' by the author).
        \end{conjecture}
        
        Suppose this were not the case, then our neat classifications of symmetric games fall to pieces as it is not intuitively clear by name that a fully symmetric game is not defined as a game which has a player $n$-transitive and strategy trivial symmetry group. 
    
        The lack of the desirable property for VNM symmetric games outlined above is a result of the definition implicitly imposing an order on each player's strategy set. We now look at where VNM symmetric games lie with respect to our other classifications of symmetric games. 
        
        \begin{proposition} 
            A game is VNM symmetric if and only if for each player permutation $\pi \in S_N$ we have $u_i(s_1, ..., s_n) = u_{\pi(i)}(s_{\pi^{-1}(1)}, ..., s_{\pi^{-1}(n)})$ for all $i \in N$ and $(s_1, ..., s_n) \in A$.
            
            \begin{proof}
                Let $\Gamma$ be a VNM symmetric game, and let $i \in N$ and $(s_1, ..., s_n) \in A$. It follows that $u_{\pi(i)}(s_1, ..., s_n) = u_i(s_{\pi(1)}, ..., s_{\pi(n)})$ for all $\pi \in S_N$. Therefore $u_{i}(s_1, ..., s_n) = u_{\pi^{-1}(i)}(s_{\pi(1)}, ..., s_{\pi(n)})$ for all $\pi \in S_N$. Subbing in $\pi = \pi^{-1}$ we have $u_i(s_1, ..., s_n) = u_{\pi(i)}(s_{\pi^{-1}(1)}, ..., s_{\pi^{-1}(n)})$ for all $\pi \in S_N$. 
                
                The converse is the same in reverse.
            \end{proof}
        \end{proposition}
    
        \begin{proposition} 
            If a game is VNM symmetric then it has a player $n$-transitive and strategy trivial symmetry group.
            
            \begin{proof}
                Let $\Gamma$ be a VNM symmetric game, $j \in N$ and $M = \{(s_j, ..., s_j): s_j \in A_j\}$. Then for each $\pi \in S_N$ we have, $u_i(s) = u_i(s_1, ..., s_n) = u_{\pi(i)}(s_{\pi^{-1}(1)}, ..., s_{\pi^{-1}(n)}) = u_{M_{\pi}.i}(M_{\pi}.s)$ for all $s \in A$.
                
                Therefore $\{M_{\pi} \in \Aut(\Gamma): \pi \in S_N\}$ is a player $n$-transitive and strategy trivial symmetry group.
            \end{proof}
        \end{proposition}
        
        \begin{proposition} 
            If a game has a player $n$-transitive and strategy trivial symmetry group then it is isomorphic to a VNM symmetric game.
            
            \begin{proof}
                If $\Gamma$ is a game with a player $n$-transitive and strategy trivial symmetry group $G$ then there exists a matching of the strategy sets $M$ such that $G = \{M_{\pi}: \pi \in S_N\}$. If we relabel the elements of each $n$-tuple in $M$ as the same then the resulting game is VNM symmetric.
            \end{proof}
        \end{proposition}
        
        This shows that VNM symmetric games have the same structure as those with a player $n$-transitive and strategy trivial symmetry group. If our above conjecture is true then VNM symmetric games have the same structure as fully standard symmetric games.
        
        A game with a player $n$-transitive and strategy trivial symmetry group has the additional property that each player is playing each of their opponents symmetrically regardless of how all other opponents play. This means these games are fair regardless of the skill for each player.

    \section{Classifying a Game}
        While the notions of symmetric, standard symmetric and fully symmetric games give us various descriptive definitions of strategic fairness, they do not provide a constructive way to determine where a particular game lies in these classifications. Below we provide an overview of some ways to classify a given game.
        
        For a game to be fully symmetric it is necessary that there exist an automorphism for every player permutation. Therefore, to show that a game is not fully symmetric we need to show there does not exist an automorphism for some player permutation, and to show that a game is fully symmetric we need to show there exists an automorphism for each player permutation in some generating set of the player permutations.
        
        Since $S_2$ is its only transitive subgroup, it follows that a two player game is symmetric if and only if it is fully symmetric. Therefore to show whether or not a two player game is symmetric we can just show whether or not it is fully symmetric.
        
        Since $<(123)>$ is the least transitive subgroup of $S_3$ when ordered by $\subseteq$, for a three player game to be symmetric it is necessary that there exist an automorphism for each 3-cycle player permutation. Therefore, to show that a three player game is not symmetric we need to show there does not exist an automorphism for some 3-cycle player permutation, and to show that a three player game is symmetric we need to show there exists an automorphism for a 3-cycle player permutation.
        
        In general it is a little more difficult to show whether or not a game is symmetric. To show that a game is symmetric we need to show there exists an automorphism for every player permutation in some generating set of some transitive subgroup of the player permutations, and to show that a game is not symmetric we must either find $\Aut(\Gamma)$ and show that it is not player transitive or show that there does not exist an automorphism for some player permutation in each generating set of every transitive subgroup of the player permutations.
        
        We note that every $n$-cycle generates a transitive subgroup of $S_n$ but not all subgroups of $S_n$ contain an $n$-cycle for $n \geq 4$. For example, $\{e, (12)\circ (34), (13)\circ (24), (14) \circ (23)\}$ is a transitive subgroup of $S_4$ which does not contain a 4-cycle. 
        
        To show that a game is standard symmetric we need to find a matching of the strategy sets such that each player permutation in a transitive subgroup $H \subseteq S_N$ is an automorphism. If there does not exist $d_i$ pure strategies such that each player plays each of their own pure strategies in exactly one strategy and the utilities are all equal for each of these strategies, then we can conclude that the game is not standard symmetric.

    \section{Constructing Symmetric Games}
        We can construct a symmetric game by picking out game bijections that generate a player transitive subgroup $G$ of the bijections and for each $g \in G$ setting $u_i(s) = u_{g.i}(g.s)$ for all $i \in N$ and $s \in A$. 
        
        \begin{example}
            If we take the following as our subgroup of bijections.
            \[ G = \{ \bigl(e ; \bigl(\begin{smallmatrix} a & b \\ a & b \end{smallmatrix}\bigr), \bigl(\begin{smallmatrix} c & d \\ c & d \end{smallmatrix}\bigr)\bigr),
                                \bigl((12) ; \bigl(\begin{smallmatrix} a & b \\ c & d \end{smallmatrix}\bigr), \bigl(\begin{smallmatrix} c & d \\ a & b \end{smallmatrix}\bigr)\bigr) \}  \]
            
            Then $\bigl((12) ; \bigl(\begin{smallmatrix} d & c \\ b & a \end{smallmatrix}\bigr), \bigl(\begin{smallmatrix} d & c \\ a & b \end{smallmatrix}\bigr)\bigr) \in G$ requires that we have,
            \begin{align*}
                u_1(a, c) &= u_2(a, c) = \alpha   &   u_1(a, d) &= u_2(b, c) = \gamma \\
                u_1(b, c) &= u_2(a, d) = \beta   &   u_1(b, d) &= u_2(b, d) = \delta
            \end{align*} 
            
            Which results in the game given below.
            \begin{center}
                \begin{game}{2}{2}
                            \> $c$  \> $d$ \\
                    $a$   \> $\alpha, \alpha $  \> $\gamma, \beta $ \\
                    $b$   \> $\beta, \gamma $  \> $\delta, \delta $
                \end{game}
            \end{center}
            
            All $2 \times 2$ standard symmetric games are isomorphic to a game of this form, hence this essentially gives us a general form for the $2 \times 2$ standard symmetric games. We can pick out a particular $2 \times 2$ standard symmetric game by assigning values to $\alpha, \beta, \gamma$ and $\delta$.
        \end{example}
    
        To construct a game which is or is not fully symmetric we just need to choose bijections such that the subgroup they generate is or is not player $n$-transitive. 
        
        To construct a standard symmetric game we can pick out a matching of the strategy sets and player permutations in a transitive subgroup of $S_n$, then the game constructed from the induced game bijections is standard symmetric. Constructing a non-standard symmetric game is not so simple, we must actually construct the game and check that it is non-standard symmetric.
    
    \section{Examples of Symmetric Games}
        Below we provide some less trivial examples in each class of symmetric games which were generated using the construction method outlined in the last section.
    
        \begin{example}
            Fully standard symmetric three player game
            \begin{center}
                \begin{game}{2}{2}[$e$]
                          \> $c$      \> $d$ \\
                    $a$   \> $1,1,1$  \> $6,2,6$ \\
                    $b$   \> $2,6,6$  \> $5,5,3$
                \end{game}
                \hspace*{10mm} 
                \begin{game}{2}{2}[$f$]
                          \> $c$     \> $d$ \\
                    $a$   \> $6,6,2$ \> $3,5,5$ \\
                    $b$   \> $5,3,5$ \> $4,4,4$
                \end{game}
            \end{center}
            
            If we take the matching $\{(a, c, e), (b,d,f)\}$ of the strategy sets then the induced bijection for each player permutation is an automorphism. The set of these automorphisms forms a player transitive and strategy trivial symmetry group making it a fully standard symmetric game.
        \end{example}

        \begin{example}
            Non-fully standard symmetric three player game.
            \begin{center}
                \begin{game}{2}{2}[$e$]
                          \> $c$      \> $d$ \\
                    $a$   \> $1,1,1$  \> $3,5,7$ \\
                    $b$   \> $5,7,3$  \> $2,4,6$
                \end{game}
                \hspace*{10mm} 
                \begin{game}{2}{2}[$f$]
                          \> $c$     \> $d$ \\
                    $a$   \> $7,3,5$ \> $6,2,4$ \\
                    $b$   \> $4,6,2$ \> $8,8,8$
                \end{game}
            \end{center}

            \begin{align*}
                \Aut(\Gamma) = \{ \bigl(e ; \bigl(\begin{smallmatrix} a & b \\ a & b \end{smallmatrix}\bigr),& \bigl(\begin{smallmatrix} c & d \\ c & d \end{smallmatrix}\bigr), \bigl(\begin{smallmatrix} e & f \\ e & f \end{smallmatrix}\bigr)\bigr),
                               \bigl((123) ; \bigl(\begin{smallmatrix} a & b \\ c & d \end{smallmatrix}\bigr), \bigl(\begin{smallmatrix} c & d \\ e & f \end{smallmatrix}\bigr), \bigl(\begin{smallmatrix} e & f \\ a & b \end{smallmatrix}\bigr)\bigr), \\
                               &\bigl((321) ; \bigl(\begin{smallmatrix} a & b \\ e & f \end{smallmatrix}\bigr), \bigl(\begin{smallmatrix} c & d \\ a & b \end{smallmatrix}\bigr), \bigl(\begin{smallmatrix} e & f \\ c & d \end{smallmatrix}\bigr)\bigr),\} 
            \end{align*}
           
            Clearly $\Aut(\Gamma)$ is player transitive and strategy trivial making $\Gamma$ standard symmetric. Additionally, since there does not exist an automorphism for every permutation of the players, $\Gamma$ is non-fully standard symmetric.
            
        \end{example}
    
        \begin{example}
            Two non-fully non-standard symmetric four player games. The second is an example of a symmetric four player game which does not have an automorphism for each 4-cycle player permutation.
            
            \begin{center}
                \begin{game}{2}{2}[$(e,g)$]
                          \> $c$      \> $d$ \\
                    $a$   \> $1,2,3,4$  \> $4,1,3,2$ \\
                    $b$   \> $2,3,4,1$  \> $5,6,7,8$
                \end{game}
                \hspace*{10mm} 
                \begin{game}{2}{2}[$(f,g)$]
                          \> $c$     \> $d$ \\
                    $a$   \> $8,5,6,7$ \> $3,4,1,2$ \\
                    $b$   \> $7,8,5,6$ \> $6,7,8,5$
                \end{game}
                \\
                \begin{game}{2}{2}[$(e,h)$]
                          \> $c$      \> $d$ \\
                    $a$   \> $6,7,8,5$  \> $7,8,5,6$ \\
                    $b$   \> $3,4,1,2$  \> $8,5,6,7$
                \end{game}
                \hspace*{10mm} 
                \begin{game}{2}{2}[$(f,h)$]
                          \> $c$     \> $d$ \\
                    $a$   \> $5,6,7,8$ \> $2,3,4,1$ \\
                    $b$   \> $4,1,2,3$ \> $1,2,3,4$
                \end{game}
                
            \end{center}
            
           \[ \Aut(\Gamma) = < \bigl((1234) ; \bigl(\begin{smallmatrix} a & b \\ d & c \end{smallmatrix}\bigr), \bigl(\begin{smallmatrix} c & d \\ e & f \end{smallmatrix}\bigr), \bigl(\begin{smallmatrix} e & f \\ g & h \end{smallmatrix}\bigr), \bigl(\begin{smallmatrix} g & h \\ a & b \end{smallmatrix}\bigr)\bigr) >  \]

            \begin{center}
                \begin{game}{2}{2}[$(e,g)$]
                          \> $c$      \> $d$ \\
                    $a$   \> $1,1,2,2$  \> $3,4,4,3$ \\
                    $b$   \> $4,3,3,4$  \> $2,2,1,1$
                \end{game}
                \hspace*{10mm} 
                \begin{game}{2}{2}[$(f,g)$]
                          \> $c$      \> $d$ \\
                    $a$   \> $4,3,3,4$  \> $2,2,1,1$ \\
                    $b$   \> $1,1,2,2$  \> $3,4,4,3$
                \end{game}
                \\
                \begin{game}{2}{2}[$(e,h)$]
                          \> $c$     \> $d$ \\
                    $a$   \> $3,4,4,3$ \> $1,1,2,2$ \\
                    $b$   \> $2,2,1,1$ \> $4,3,3,4$
                \end{game}
                \hspace*{10mm} 
                \begin{game}{2}{2}[$(f,h)$]
                          \> $c$     \> $d$ \\
                    $a$   \> $2,2,1,1$ \> $4,3,3,4$ \\
                    $b$   \> $3,4,4,3$ \> $1,1,2,2$
                \end{game}
                
            \end{center}
            
            {\footnotesize 
            \begin{align*}
                \Aut(\Gamma) &= < \bigl((12)(34) ; \bigl(\begin{smallmatrix} a & b \\ d & c \end{smallmatrix}\bigr), \bigl(\begin{smallmatrix} c & d \\ a & b \end{smallmatrix}\bigr), \bigl(\begin{smallmatrix} e & f \\ h & g \end{smallmatrix}\bigr), \bigl(\begin{smallmatrix} g & h \\ e & f \end{smallmatrix}\bigr)\bigr), \\
                               &\bigl((13)(24) ; \bigl(\begin{smallmatrix} a & b \\ f & e \end{smallmatrix}\bigr), \bigl(\begin{smallmatrix} c & d \\ h & g \end{smallmatrix}\bigr), \bigl(\begin{smallmatrix} e & f \\ a & b \end{smallmatrix}\bigr), \bigl(\begin{smallmatrix} g & h \\ c & d \end{smallmatrix}\bigr)\bigr), 
                               \bigl((14)(23) ; \bigl(\begin{smallmatrix} a & b \\ h & g \end{smallmatrix}\bigr), \bigl(\begin{smallmatrix} c & d \\ f & e \end{smallmatrix}\bigr), \bigl(\begin{smallmatrix} e & f \\ c & d \end{smallmatrix}\bigr), \bigl(\begin{smallmatrix} g & h \\ a & b \end{smallmatrix}\bigr)\bigr) >
           \end{align*}
           }
           
           Clearly the game bijections generate games which are symmetric and not fully symmetric. Additionally, since there exist no strategy profiles such that the players' utility values are equal, the games are not standard symmetric. Hence, both games are non-fully non-standard symmetric.
        \end{example}
            
        \begin{example}
            Fully non-standard symmetric four player game.
            
            \begin{center}
                \begin{game}{2}{2}[$(e,g)$]
                          \> $c$      \> $d$ \\
                    $a$   \> $2,1,1,1$  \> $1,1,2,1$ \\
                    $b$   \> $2,1,1,1$  \> $1,1,1,2$
                \end{game}
                \hspace*{10mm} 
                \begin{game}{2}{2}[$(f,g)$]
                          \> $c$      \> $d$ \\
                    $a$   \> $1,1,1,2$  \> $1,1,2,1$ \\
                    $b$   \> $1,2,1,1$  \> $1,2,1,1$
                \end{game}
                \\
                \begin{game}{2}{2}[$(e,h)$]
                          \> $c$     \> $d$ \\
                    $a$   \> $1,2,1,1$ \> $1,2,1,1$ \\
                    $b$   \> $1,1,2,1$ \> $1,1,1,2$
                \end{game}
                \hspace*{10mm} 
                \begin{game}{2}{2}[$(f,h)$]
                          \> $c$     \> $d$ \\
                    $a$   \> $1,1,1,2$ \> $2,1,1,1$ \\
                    $b$   \> $1,1,2,1$ \> $2,1,1,1$
                \end{game}
                
            \end{center}
            
            {\small
                \[ \Aut(\Gamma) = < \bigl((1234) ; \bigl(\begin{smallmatrix} a & b \\ c & d \end{smallmatrix}\bigr), \bigl(\begin{smallmatrix} c & d \\ e & f \end{smallmatrix}\bigr), \bigl(\begin{smallmatrix} e & f \\ h & g \end{smallmatrix}\bigr), \bigl(\begin{smallmatrix} g & h \\ a & b \end{smallmatrix}\bigr)\bigr),
                               \bigl((12) ; \bigl(\begin{smallmatrix} a & b \\ c & d \end{smallmatrix}\bigr), \bigl(\begin{smallmatrix} c & d \\ a & b \end{smallmatrix}\bigr), \bigl(\begin{smallmatrix} e & f \\ e & f \end{smallmatrix}\bigr), \bigl(\begin{smallmatrix} g & h \\ h & g \end{smallmatrix}\bigr)\bigr) > \]
            }
            
            Clearly the game bijections generate a fully symmetric game, and since there exist no strategy profiles such that the players' utility values are equal the game is not standard symmetric.
        \end{example}
        
    \section{Incorrect Symmetric Game Classification}      
        Dasgupta and Maskin \cite{DMaskin} gave an incorrect definition of a symmetric game which we give below.
        
        \begin{definition} 
            A game $\Gamma$ is \textit{DM symmetric} if $A_i = A_j$ for all $i, j \in N$, and for each player permutation $\pi \in S_N$,  $u_i(s_1, ..., s_n) = u_{\pi(i)}(s_{\pi(1)}, ..., s_{\pi(n)})$ for all $i \in N$ and $(s_1, ..., s_n) \in A$.
        \end{definition}
        
        The problem with this definition is that it does not correctly permute the players and their strategies, for example for each player $i$, the right hand side does not have player $\pi(i)$ playing the strategy that player $i$ played.
    
        \begin{proposition} 
            A game is VNM symmetric if and only if for each player transposition $\pi \in S_N$, $u_i(s_1, ..., s_n) = u_{\pi(i)}(s_{\pi(1)}, ..., s_{\pi(n)})$ for all $i \in N$ and $(s_1, ..., s_n) \in A$.
            
            \begin{proof}
                Let $\Gamma$ be a VNM symmetric game and let $\pi \in S_N$ be a transposition. Then we have $u_i(s_1, ..., s_n) = u_{\pi(i)}(s_{\pi^{-1}(1)}, ..., s_{\pi^{-1}(n)})$ for all $i \in N$ and $(s_1, ..., s_n) \in A$. Therefore since $\pi = \pi^{-1}$ we have $u_i(s_1, ..., s_n) = u_{\pi(i)}(s_{\pi(1)}, ..., s_{\pi(n)})$ for all $i \in N$ and $(s_1, ..., s_n) \in A$.
                
                The converse is much the same except we note that the condition for a DM symmetric game is met for all player transpositions since it is met for all player permutations.
            \end{proof}
        \end{proposition}
            
        This shows that if a game is DM symmetric then it is VNM symmetric. Furthermore, since each element of $S_2$ is its own inverse, it follows that the two notions are equivalent for two player games. 
        
        Below we give an example of a three player VNM symmetric game which is not DM symmetric to show that the classifications are not equivalent for games with at least three players.
    
        \begin{example}
            Three Player Symmetric Game
            \begin{center}
                \begin{game}{2}{2}[$a$]
                          \> $a$      \> $b$ \\
                    $a$   \> $1,1,1$  \> $2,3,2$ \\
                    $b$   \> $3,2,2$  \> $4,4,5$
                \end{game}
                \hspace*{10mm} 
                \begin{game}{2}{2}[$b$]
                          \> $a$     \> $b$ \\
                    $a$   \> $2,2,3$ \> $5,4,4$ \\
                    $b$   \> $4,5,4$ \> $6,6,6$
                \end{game}
            \end{center}
            
            If we take $\pi = \bigl(\begin{smallmatrix} 1 & 2 & 3 \\ 2 & 3 & 1 \end{smallmatrix}\bigr) \in S_3$ and $s = (b, a, a) \in A$, we see that
            \begin{align*}
                3 = u_1(b, a, a) &= u_1(s_1, s_2, s_3) \\
                & \neq u_{\pi(1)}(s_{\pi(1)}, s_{\pi(2)}, s_{\pi(3)}) = u_2(s_2, s_3, s_1) = u_2(a, a, b) = 2
            \end{align*} 
            
            It should be fairly obvious that if we are mapping player 1 to player 2 and player 1 is playing $b$ then we want the mapped strategy to have player 2 playing $b$, but this is not what is required for a DM symmetric game.
        \end{example}
            
        We could restrict the condition for a DM symmetric game to just player transpositions, however transpositions are not closed under composition so the composition of automorphisms would not necessarily meet this condition, making it fairly undesirable. 
        
    \section{Equilibria in Symmetric Games}
        We now present two famous theorems without proof relating to the existence of Nash equilibria in symmetric games. Our first result from Nash \cite{NashNCG} gives the existence of a symmetric Nash equilibrium.
        
        \begin{theorem}
            Every symmetric game has a symmetric Nash equilibrium.
        \end{theorem}
        
        Our second result by Cheng et al. \cite{CRVWSym} gives the existence of a pure strategy Nash equilibrium for fully standard symmetric $2 \times ... \times 2$ games.
        
        \begin{theorem} 
            Every fully standard symmetric $2\times ... \times 2$ game has a pure strategy Nash equilibrium.
        \end{theorem}
        
        Cheng \textit{et al.} \cite{CRVWSym} noted that rock, paper, scissors is an example of a fully standard symmetric $3 \times 3$ game with no pure strategy Nash equilibria, and indirectly that matching pennies is an example of a fully non-standard symmetric $2\times 2$ game with no pure strategy Nash equilibria. 
        
        \begin{example}
            Rock, paper, scissors.
            
            \begin{center}
                \begin{game}{3}{3}
                          \> $R$      \> $P$      \> $S$\\
                    $R$   \> $0,0$  \> $0,1$  \> $1,0$ \\
                    $P$   \> $1,0$  \> $0,0$  \> $0,1$ \\
                    $S$   \> $0,1$  \> $1,0$  \> $0,0$
                \end{game}
            \end{center}
        \end{example}
        
        However they only considered fully standard symmetric games as symmetric, consequently not considering whether the result might hold for non-fully standard symmetric games. We provide an example below to show that this is not the case.
        
        \begin{example}
            A non-fully standard symmetric $2 \times 2 \times 2$ game with no pure strategy Nash equilibria.
            
            \begin{center}
                \begin{game}{2}{2}[$e$]
                          \> $c$      \> $d$ \\
                    $a$   \> $2,2,2$  \> $6,4,5$ \\
                    $b$   \> $4,5,6$  \> $8,1,7$
                \end{game}
                \hspace*{10mm} 
                \begin{game}{2}{2}[$f$]
                          \> $c$     \> $d$ \\
                    $a$   \> $5,6,4$ \> $7,8,1$ \\
                    $b$   \> $1,7,8$ \> $3,3,3$
                \end{game}
            \end{center}
        \end{example}

    \bibliographystyle{plain}
    \bibliography{references}
\end{document}